\documentclass[10pt]{article} \usepackage{amsfonts} \usepackage{amssymb}\usepackage{amsmath}\usepackage{amsfonts}
\usepackage[russian]{babel}

\makeatletter
\renewcommand{\@makefnmark}{}
\makeatother
\begin{document}
\baselineskip=10pt
\pagestyle{plain}
{\Large
\newcommand{\mun}{\mu_{n,j}}
\newcommand{\wnj}{W_{N,j}(\mu)}

\newcommand{\ntm}{|2n-\theta-\mu|}
\newcommand{\pwp}{PW_\pi^-}
\newcommand{\cpm}{\cos\pi\mu}
\newcommand{\spm}{\frac{\sin\pi\mu}{\mu}}
\newcommand{\tet}{(-1)^{\theta+1}}
\newcommand{\dm}{\Delta(\mu)}
\newcommand{\dl}{\Delta(\lambda)}

\newcommand{\dnm}{\Delta_N(\mu)}
\newcommand{\sni}{\sum_{n=N+1}^\infty}
\newcommand{\muk}{\sqrt{\mu^2+q_0}}
\newcommand{\agt}{\alpha, \gamma, \theta,}
\newcommand{\muq}{\sqrt{\mu^2-q_0}}
\newcommand{\smn}{\sum_{n=1}^\infty}
\newcommand{\lop}{{L_2(0,2\pi)}}
\newcommand{\tdm}{\tilde\Delta_+(\mu_n)}
\newcommand{\dmn}{\Delta_+(\mu_n)}
\newcommand{\dd}{D_+(\mu_n)}
\newcommand{\ddd}{\tilde D_+(\mu_n)}
\newcommand{\sn}{\sum_{n=1}^\infty}
\newcommand{\pnn}{\prod_{n=1}^\infty}
\newcommand{\aln}{\alpha_n(\mu)}
\newcommand{\emp}{e^{\pi|Im\mu|}}
\newcommand{\ppp}{\prod_{p=p_0}^\infty}
\newcommand{\mm}{m-1/2}
\newcommand{\pps}{\sum_{p=p_0}^\infty}
\newcommand{\mln}{m(\lambda_n)}
\newcommand{\mnk}{\mu_{n_k}}
\newcommand{\lnk}{\lambda_{n_k}}
\newcommand{\lqr}{\langle q \rangle}
\newcommand{\mlnk}{m(\lambda_{n_k})}
\newcommand{\emx}{e^{i\mu x}}
\newcommand{\emxx}{e^{-i\mu x}}

\newcommand{\om}{O(\frac{1}{\mu})}
\newcommand{\lp}{L_2(0,\pi)}
\newcommand{\lpp}{L_2(\Omega)}

\newcommand{\sxm}{s(x,\mu)}
\newcommand{\skm}{s(\xi,\mu)}
\newcommand{\cxm}{c(x,\mu)}
\newcommand{\ckm}{c(\xi,\mu)}
\newcommand{\kp}{K_{\pi/2}}

\newcommand{\No}{\textnumero}

\newcommand{\elx}{e^{i\lambda x}}
\newcommand{\elxx}{e^{-i\lambda x}}

\medskip
\medskip
\medskip

\centerline {\bf On Convergence of Spectral Expansions of Dirac Operators}
\centerline {\bf   with Regular Boundary Conditions}
\medskip
\medskip
\centerline {\bf Alexander Makin}
\medskip
\medskip
\medskip

\medskip
\medskip
\begin{quote}{\normalsize

  Spectral problem for the Dirac operator with regular but not strongly regular boundary conditions and complex-valued
  potential summable over a finite interval is considered.  The purpose of this paper is to find conditions under which the root function system forms a usual Riesz basis rather than a Riesz basis with parentheses.}
   \end{quote}

\medskip
\medskip	
\medskip

We study spectral problem for the Dirac operator, generated be the differential expression
 $$
 L\mathbf{y}=B\mathbf{y}'+V\mathbf{y},
 $$
 where
$$
 B=\begin{pmatrix}
 -i&0\\
 0&i
 \end{pmatrix},\quad V=\begin{pmatrix}
 0&P(x)\\
 Q(x)&0
 \end{pmatrix},
$$
functions $P(x), Q(x)\in L_1(0,\pi)$ and two-point boundary conditions
$$
U(\mathbf{y})=C\mathbf{y}(0)+D\mathbf{y}(\pi)=0,
$$
where
$$
C=\begin{pmatrix}
 a_{11}&a_{12}\\
 a_{21}&a_{22}
 \end{pmatrix},\quad D=\begin{pmatrix}
a_{13} &a_{14}\\
a_{23}&a_{24} \end{pmatrix},
$$
the coefficients $a_{ij}$ are arbitrary complex numbers,
and rows of  the matrix
$$
A=\begin{pmatrix}
 a_{11}&a_{12}&a_{13} &a_{14}\\
 a_{21}&a_{22}&a_{23} &a_{24}
\end{pmatrix}
$$
are linearly independent. Denote by
$$
E(x,\lambda)=\begin{pmatrix}
e_{11}(x,\lambda)&e_{12}(x,\lambda)\\
e_{21}(x,\lambda)&e_{22}(x,\lambda)
\end{pmatrix}
$$
the matrix of the fundamental solution system  of the equation
$L\mathbf{y}=\lambda\mathbf{y}$
with boundary condition
$
E(0,\lambda)=I
$, where
$I$ is the unit matrix,
and by
$$
E_0(x,\lambda)=\begin{pmatrix}
e_{11}^0(x,\lambda)&e_{12}^0(x,\lambda)\\
e_{21}^0(x,\lambda)&e_{22}^0(x,\lambda)
\end{pmatrix}
$$
the fundamental solution system  of the equation
$L_0\mathbf{y}=\lambda\mathbf{y}
$
with boundary condition
$
E_0(0,\lambda)=I$, where
$L_0\mathbf{y}=B\mathbf{y}'$.

Obviously,  $e_{11}^0(x,\lambda)=e^{i\lambda x}$, $e_{22}^0(x,\lambda)=e^{-i\lambda x}$, $e_{12}^0(x,\lambda)=e_{21}^0(x,\lambda)=0$.
Denote by $A_{ij}$ the determinant composed of the $i$th and $j$th columns  of the matrix $A$. It is known [4], that the characteristic determinant
$\Delta(\lambda)$ of problem
$$
L\mathbf{y}=\lambda\mathbf{y},\quad U(\mathbf{y})=0 \eqno(1)
$$
can be reduced to the form
$$
\begin{array}{c}
\Delta(\lambda)=A_{12}+A_{34}+A_{32}e_{11}(\pi,\lambda)+A_{14}e_{22}(\pi,\lambda)+A_{13}e_{12}(\pi,\lambda)+A_{42}e_{21}(\pi,\lambda)=\\
=\Delta_0(\lambda)+\int_0^\pi r_1(t)e^{-i\lambda t}dt+\int_0^\pi r_2(t)e^{i\lambda t}dt,
\end{array}\eqno(2)
$$
where the function

$$
\Delta_0(\lambda)=(A_{12}+A_{34})-A_{23}e^{i\pi\lambda}+A_{14}e^{-i\pi\lambda}\eqno(3)
$$
is the characteristic determinant of problem
$$
L_0\mathbf{y}=\lambda\mathbf{y},\quad U(\mathbf{y})=0, \eqno(4)
$$
functions $r_j(t)\in L_1(0,\pi)$, $j=1,2$.
Denote
$\Omega=(0,\pi)\times(0,\pi)$,
$\mathbb{H}=L_2(0,\pi)\oplus L_2(0,\pi)$,
$$
\mathcal{E}_{j1}(a,x,\lambda)=e_{j1}(x,\lambda)e_{22}(a,\lambda)-e_{j2}(x,\lambda)e_{21}(a,\lambda),
$$
$$
\mathcal{E}_{j2}(a,x,\lambda)=e_{j2}(x,\lambda)e_{11}(a,\lambda)-e_{j1}(x,\lambda)e_{12}(a,\lambda),
$$
$$
\mathcal{E}_{j1}^0(a,x,\lambda)=e_{j1}^0(x,\lambda)e_{22}^0(a,\lambda)-e_{j2}^0(x,\lambda)e_{21}^0(a,\lambda),
$$
$$
\mathcal{E}_{j2}^0(a,x,\lambda)=e_{j2}^0(x,\lambda)e_{11}^0(a,\lambda)-e_{j1}^0(x,\lambda)e_{12}^0(a,\lambda)
$$
$(j=1,2)$.
The Green function of problem (1) admits [2] the following representation
$$
\begin{array}{c}
G(t,x,\lambda)=\frac{i}{\Delta(\lambda)}\Bigl(A_{12}\begin{pmatrix}
\mathcal{E}_{11}(t,x,\lambda)&-\mathcal{E}_{12}(t,x,\lambda)\\
\mathcal{E}_{21}(t,x,\lambda)&-\mathcal{E}_{22}(t,x,\lambda)
\end{pmatrix}+\\
+\begin{pmatrix}
\mathcal{E}_{11}(\pi,x,\lambda)&-\mathcal{E}_{12}(\pi,x,\lambda)\\
\mathcal{E}_{21}(\pi,x,\lambda)&-\mathcal{E}_{22}(\pi,x,\lambda)
\end{pmatrix}\begin{pmatrix}
A_{14}&A_{24}\\
A_{13}&A_{23}
\end{pmatrix}\begin{pmatrix}
e_{22}(t,\lambda)&e_{12}(t,\lambda)\\
-e_{21}(t,\lambda)&-e_{11}(t,\lambda)
\end{pmatrix}-\\
-\Delta(\lambda)\chi_{t>x}(t,x)\begin{pmatrix}
\mathcal{E}_{11}(t,x,\lambda)&-\mathcal{E}_{12}(t,x,\lambda)\\
\mathcal{E}_{21}(t,x,\lambda)&-\mathcal{E}_{22}(t,x,\lambda)
\end{pmatrix}\Bigr),
\end{array}\eqno(5)
$$
where $\chi_{t>x}$ is the characteristic function of the triangle  $t>x$.

Multiplying the matrixes in right-hand side of (5), we obtain
$$
G(t,x,\lambda)=\frac{iH(t,x,\lambda)}{\Delta(\lambda)}-i\chi_{t>x}(t,x)\begin{pmatrix}
\mathcal{E}_{11}(t,x,\lambda)&-\mathcal{E}_{12}(t,x,\lambda)\\
\mathcal{E}_{21}(t,x,\lambda)&-\mathcal{E}_{22}(t,x,\lambda)
\end{pmatrix},\eqno(6)
$$
where there matrix
$
H(t,x,\lambda)=||h_{jk}(t,x,\lambda)||
$
has the elements
$$
\begin{array}{c}
h_{11}(t,x,\lambda)=A_{12}\mathcal{E}_{11}(t,x,\lambda)+[A_{14}\mathcal{E}_{11}(\pi,x,\lambda)-A_{13}\mathcal{E}_{12}(\pi,x,\lambda)]e_{22}(t,\lambda)-\\
-[A_{24}\mathcal{E}_{11}(\pi,x,\lambda)-A_{23}\mathcal{E}_{12}(\pi,x,\lambda)]e_{21}(t,\lambda),
\end{array}
$$
$$
\begin{array}{c}
h_{12}(t,x,\lambda)=-A_{12}\mathcal{E}_{12}(t,x,\lambda)+[A_{14}\mathcal{E}_{11}(\pi,x,\lambda)-A_{13}\mathcal{E}_{12}(\pi,x,\lambda)]e_{12}(t,\lambda)-\\
-[A_{24}\mathcal{E}_{11}(\pi,x,\lambda)-A_{23}\mathcal{E}_{12}(\pi,x,\lambda)]e_{11}(t,\lambda),
\end{array}
$$
$$
\begin{array}{c}
h_{21}(t,x,\lambda)=A_{12}\mathcal{E}_{21}(t,x,\lambda)+[A_{14}\mathcal{E}_{21}(\pi,x,\lambda)-A_{13}\mathcal{E}_{22}(\pi,x,\lambda)]e_{22}(t,\lambda)-\\
-[A_{24}\mathcal{E}_{21}(\pi,x,\lambda)-A_{23}\mathcal{E}_{22}(\pi,x,\lambda)]e_{21}(t,\lambda),
\end{array}
$$
$$
\begin{array}{c}
h_{22}(t,x,\lambda)=-A_{12}\mathcal{E}_{22}(t,x,\lambda)+[A_{14}\mathcal{E}_{21}(\pi,x,\lambda)-A_{13}\mathcal{E}_{22}(\pi,x,\lambda)]e_{12}(t,\lambda)-\\
-[A_{24}\mathcal{E}_{21}(\pi,x,\lambda)-A_{23}\mathcal{E}_{22}(\pi,x,\lambda)]e_{11}(t,\lambda).
\end{array}
$$

In the strip $\Pi:|Im\lambda|<C$  the estimates
$$
e_{jk}(x,\lambda)=e_{jk}^0(x,\lambda)+o(1),\quad
\mathcal{E}_{jk}(t,x,\lambda)=\mathcal{E}_{jk}^0(t,x,\lambda)+o(1),\eqno(7)
$$
are valid [2] uniformly by $t,x$, $j=1,2,\, k=1,2$,

It follows from [2] that the Green function $G_0(t,x,\lambda)$ of problem (4) can be represented in the form
$$
G_0(t,x,\lambda)=\frac{iH_0(t,x,\lambda)}{\Delta_0(\lambda)}-i\chi_{t>x}(t,x)\begin{pmatrix}
e^{i\lambda(x-t)}&0\\
0&-e^{i\lambda(t-x)}
\end{pmatrix},
$$
where
$$
\begin{array}{c}
H_0(t,x,\lambda)=||h_{jk}^0(t,x,\lambda)||=\\
=\begin{pmatrix}
A_{12}e^{i\lambda(x-t)}+A_{14}e^{i\lambda(x-\pi-t)}&-A_{24}e^{i\lambda(x-\pi+t)}\\
-A_{13}e^{i\lambda(\pi-x-t)}&-A_{12}e^{i\lambda(t-x)}+A_{23}e^{i\lambda(\pi-x+t)}
\end{pmatrix}=\\
=\begin{pmatrix}
e^{i\lambda(x-t)}(A_{12}+A_{14}e^{-i\pi\lambda})&-A_{24}e^{i\lambda(x-\pi+t)}\\
-A_{13}e^{i\lambda(\pi-x-t)}&e^{i\lambda(t-x)}(-A_{12}+A_{23}e^{i\pi\lambda})
\end{pmatrix}.
\end{array}\eqno(8)
$$
It follows from (7)  that the relations
$$
h_{jk}(t,x,\lambda)=h_{jk}^0(t,x,\lambda)+o(1)\eqno(9)
$$
hold  in the strip $\Pi$.

Boundary conditions $U(\mathbf{y})=0$ are called regular if
$$
A_{14}A_{23}\ne0,\eqno(10)
$$
and ones are called strongly regular if besides (10) the following condition
$$
(A_{12}+A_{34})^2+4A_{14}A_{23}\ne0
$$
takes place.
It is well known, that the eigenvalues $\lambda_n^0$ of problem (4) with regular boundary conditions forms two series
$\lambda_n^0=\lambda_{k,1}^0=-\frac{i}{\pi}\ln z_1+2k$ if $n=2k$, and $\lambda_n^0=\lambda_{k,2}^0=-\frac{i}{\pi}\ln z_2+2k$ if $n=2k+1$, $k\in\mathbb{Z}$, where
$z_1$ and $z_2$ are the roots of square equation
$$
A_{23}z^2-(A_{12}+A_{34})z -A_{14}=0,
$$
$-\pi<Im z_j\le\pi$, $j=1,2$.
Also, it is known [2] that the eigenvalues $\lambda_n$ of problem (1) with regular boundary conditions satisfy
the asymptotic relation
$$
\lambda_n=\lambda_n^0+o(1),\eqno(11)
$$
where $\lambda_n^0$ are the eigenvalues of corresponding problem (4).

It was established in [1-4] that the root function system of problem (1) with strongly regular boundary conditions forms a Riesz basis in $\mathbb{H}$
and a Riesz basis with parentheses in $\mathbb{H}$ in the case of regular but not strongly regular boundary conditions.
Further, we consider problem (1) only with regular but not strongly regular boundary conditions, i.e., we suppose that
$$
(A_{12}+A_{34})^2+4A_{14}A_{23}=0.\eqno(12)
$$
One can readily see that in this case
$$
z_1=z_2=\frac{A_{12}+A_{34}}{2A_{23}}=z,
$$
and the spectrum consists of pairwise close eigenvalues
$\lambda_n=\lambda_{k,1}=-\frac{i}{\pi}\ln z+2k+ \varepsilon_{k,1}$ if $n=2k$, and $\lambda_n=\lambda_{k,2}=-\frac{i}{\pi}\ln z+2k+ \varepsilon_{k,2}$ if $n=2k+1$, $k\in\mathbb{Z}$, $\varepsilon_{k,j}\to0$ if $k\to\pm\infty$. If for all $k$ such that $|k|>k_0$, $\lambda_{k,1}=\lambda_{k,2}$, then the spectrum is called asymptotically multiple.

Following [5], we call regular but not strongly regular boundary conditions determined by
the matrix $A$ periodic-type if
$A_{13}=A_{24}=0,\quad A_{12}=A_{34}$.
Periodic-type boundary conditions are equivalent to the conditions given by the matrix
$$
\begin{pmatrix}
 1&0&a&0\\
0 &a&0&1
\end{pmatrix},\eqno(13)
$$
where $a\ne0$. If $a=-1$, then boundary conditions (13) are periodic, and if $a=1$, ones are antiperiodic.
It is known [5], that for problem (4) with periodic-type conditions all root subspaces consist of two eigenfunctions,
and for all remaining conditions they consist of one eigenfunction and one associated function. In all cases the root function system of indicated
problem forms a Riesz basis  in $\mathbb{H}$.

{\bf Theorem 1.} {\it If the boundary conditions are not periodic-type conditions, then the system of eigen- and associated functions of problem (1)  forms a Riesz basis in  $\mathbb{H}$ if and only if the spectrum is asymptotically multiple.}

Proof. The root function system of problem (1) is a Riesz basis of subspaces, corresponding to closed eigenvalues $\lambda_{k,1}$ and $\lambda_{k,2}$. Let the spectrum of problem (1) be asymptotically multiple.  Choosing in every two-dimensional root subspaces, corresponding to multiple eigenvalues, an orthonormal basis, we find [6, p. 414] that  the sequence of root functions, which is the union of all eigenfunctions, corresponding to simple eigenvalues and all indicated  orthonormal bases,  forms a Riesz basis in $\mathbb{H}$.

By the condition of Theorem 1,
$$
|A_{13}|+|A_{24}|+|A_{12}-A_{34}|>0.\eqno(14)
$$
If $A_{13}\ne0$, then it follows from (8) that
$$
||h_{21}^0(t,x,\lambda_n)||_{\lpp}>c_1>0,
$$
which, together with (9) implies
$$
||h_{21}(t,x,\lambda_n)||_{\lpp}>c_2>0.
$$
If $A_{24}\ne0$, the it follows from (8) that
 $$
 ||h_{12}^0(t,x,\lambda_n)||_{\lpp}>c_3>0.
 $$
 This and (9) imply
 $$
||h_{12}(t,x,\lambda_n)||_{\lpp}>c_4>0.
$$
From the equation
$\Delta_0(\lambda)=0$ and
relations (10), (12)
we obtain
$$
e^{i\pi\lambda_n^0}=\frac{A_{12}+A_{34}}{2A_{23}},\quad e^{-i\pi\lambda_n^0}=-\frac{A_{12}+A_{34}}{2A_{14}},
$$
hence,
$$
A_{12}+A_{14}e^{-i\pi\lambda_n^0}=(A_{12}-A_{34})/2,\quad
-A_{12}+A_{23}e^{i\pi\lambda_n^0}=(A_{34}-A_{12})/2.
$$
Therefore, if $|A_{12}-A_{34}|>0$, then it follows from (8) that
$$
||h_{11}^0(t,x,\lambda_n^0)||_{\lpp}>c_5>0,\quad||h_{22}^0(t,x,\lambda_n^0)||_{\lpp}>c_6>0.
$$
which, together with (11) implies
$$
||h_{11}^0(t,x,\lambda_n)||_{\lpp}>c_7>0,\quad||h_{22}^0(t,x,\lambda_n)||_{\lpp}>c_8>0.
$$
It follows from two last inequalities and (9) that
$$
||h_{11}(t,x,\lambda_n)||_{\lpp}>c_9>0,\quad||h_{22}(t,x,\lambda_n)||_{\lpp}>c_{10}>0.
$$
Thus, under condition (14) there exist $j,k$, such that  for all sufficiently large $|n|$
$$
||h_{jk}(t,x,\lambda_n)||_{\lpp}>c_0>0.\eqno(15)
$$

Let us estimate the function $\Delta'(\lambda)$. It follows from (2), (3) and  the Riemann lemma that in the strip $\Pi$
$$
\Delta'(\lambda_n)=-i\pi(A_{23}e^{i\pi\lambda_n}+A_{14}e^{-i\pi\lambda_n})+o(1)=-i\pi(A_{23}e^{i\pi\lambda_n^0}+A_{14}e^{-i\pi\lambda_n^0})+o(1)=o(1).
\eqno(16)
$$
Let $\lambda_n$ be a simple eigenvalue of problem (1),  let $\mathbf{y}_n(x)=(y_n^{[1]}(x),y_n^{[2]}(x))$ be the corresponding eigenfunction, and let $\mathbf{z}_n(x)=(z_n^{[1]}(x),z_n^{[2]}(x))$ be the function from the biorthogonally conjugate  system that form a pair with $\mathbf{y}_n(x)$. Then [7], the principal part of Green function  $G(t,x,\lambda)$ in the  neighborhood of
$\lambda_n$ has the form
$$
\frac{\mathbf{y}_n(x)\overline{\mathbf{z}_n(t)}}{\lambda-\lambda_n}.
$$
This and (6) imply
$$
\mathbf{y}_n(x)\overline{\mathbf{z}_n(t)}=\frac{iH(t,x,\lambda_n)}{\Delta'(\lambda_n)}.\eqno(17)
$$
If the spectrum of problem (1) is not asymptotically multiple, then there exists an infinite subsequence $\lambda_{n_m}$ $(m=1,2,\ldots)$ of simple eigenvalues.
It follows from (15-17) that
$$
\sum_{j=1,k=1}^2||y_{n_m}^{[j]}z_{n_m}^{[k]}||_{\lpp}^2=||\mathbf{y}_{n_m}||_\mathbb{H}^2||\mathbf{z}_{n_m}||_\mathbb{H}^2\to\infty
$$
as $m\to\infty$.
This, together with a theorem of resonance type [8, p. 104] implies existence of a function $\mathbf{f}\in\mathbb{H}$ such that
$$
|\langle\mathbf{f},\mathbf{z}_{n_m}\rangle|||\mathbf{y}_{n_m}||_\mathbb{H}^2\to\infty
$$
as $m\to\infty$.
Theorem 1 is proved.

Remark. For the periodic-type conditions
$h_{jk}^0(t,x,\lambda_n^0)\equiv0$, where $j=1,2,\,k=1,2$.

The case of periodic and antiperiodic boundary conditions was investigated in [5, 9-15] if the functions $P(x), Q(x)\in L_2(0,\pi)$. In particular, an example of potential matrix $V(x)$, providing  that the corresponding eigenfunction expansion diverges in $\mathbb{H}$ was constructed in [9, Th. 71].

Consider problem (1) with periodic-type boundary conditions. It has the form

$$
\begin{array}{c}
-iy_1'+P(x)y_2=\lambda y_1,\\
iy_2'+Q(x)y_1=\lambda y_2,\\
y_1(0)+ay_1(\pi)=0,\quad ay_2(0)+y_2(\pi)=0.
\end{array}\eqno(18)
$$
Let $a=re^{i\varphi}, -\pi<\varphi\le\pi$.
Denote $\tau_0=\frac{\varphi+\pi}{\pi}+\frac{i\ln r}{\pi}$. Then problem (18) has two series of the eigenvalues
$$
\lambda_{n,j}=\tau_0+2n+\varepsilon_{n,j},\eqno(19)
$$
where $n\in\mathbb{Z}$, $\varepsilon_{n,j}\to0$ as $|n|\to\infty$, $j=1,2$. This, in particular, implies that
if $|n|>n_0$, then the numbers $\lambda_{n,j}$ lie inside the disks of radius $1/2$ centered at the points $\tau_0+2n$.
The conjugate problem has the form
$$
\begin{array}{c}
-iz_1'+\bar Q(x)z_2=\lambda z_1,\\
iz_2'+\bar P(x)z_1=\lambda z_2,\\
\bar az_1(0)+z_1(\pi)=0,\quad z_2(0)+\bar a z_2(\pi)=0.
\end{array}\eqno(18')
$$

Let $\lambda_{n,j}$ be a simple eigenvalue of problem (18),  and let $\mathbf{y}_n(x)=(y_{n,j}^{[1]}(x),y_{n,j}^{[2]}(x))$ be the corresponding eigenfunction.
Let $\mathbf{z}_{n,j}(x)=(z_{n,j}^{[1]}(x),z_{n,j}^{[2]}(x))$ be the function from biorthogonally conjugate system that forms a pair with $\mathbf{y}_{n,j}(x)$. Since the function $(\bar y_{n,j}^{[2]}(x),\bar y_{n,j}^{[1]}(x))$ is an eigenfunction of problem $(18')$, corresponding to the eigenvalue $\bar\lambda_{n,j}$, then  $\mathbf{z}_{n,j}(x)=\gamma_{n,j}(\bar y_{n,j}^{[2]}(x),\bar y_{n,j}^{[1]}(x))$.
Since $\langle\mathbf{y}_{n,j},\mathbf{z}_{n,j}\rangle=1$, then
$$
\bar\gamma_{n,j}\int_0^\pi y_{n,j}^{[1]}(x),y_{n,j}^{[2]}(x)dx=1/2.\eqno(20)
$$

It is easily seen that
$$
 \mathbf{y}_{n,j}(x)=\alpha_{n,j}\begin{pmatrix}
e_{11}(x,\lambda_{n,j})\\
e_{21}(x,\lambda_{n,j})
\end{pmatrix}+
\beta_{n,j}\begin{pmatrix}
e_{12}(x,\lambda_{n,j})\\
e_{22}(x,\lambda_{n,j})
\end{pmatrix}.
$$
It follows from asymptotic relations (7) for the functions $e_{jk}(x,\lambda)$
that
$$
y_{n,j}^{[1]}(x)=\alpha_{n,j}(e^{i\lambda_n x}+o(1))+\beta_{n,j}o(1),\quad y_{n,j}^{[2]}(x)=\alpha_{n,j}o(1)+\beta_{n,j}(e^{-i\lambda_n x}+o(1)).\eqno(21)
$$

In what follows, we assume without loss of generality that the root functions corresponding  to a multiple eigenvalue form an orthonormal basis in the corresponding root subspace. Let $T$ denote the set of those numbers $n$ for which $\lambda_{n,1}\ne\lambda_{n,2}$. If $T$ is finite, then, according to [6, p. 414] the system of root functions is a Riesz basis in $\mathbb{H}$.
Further, we suppose that $T$ is infinite, $n\in T$.

{\bf Lemma 1.} {\it The system of eigen- and associated function of problem (18) is a Riesz basis if and only if there exist constants $C_1, C_2>0$ such that for all sufficiently large $|n|$
$$
C_1<\biggl|\frac{\alpha_{n,1}}{\beta_{n,1}}\biggl|<C_2.\eqno(22)
$$
}
Proof. We can assume without loss of generality that
$$
|\alpha_{n,j}|^2+|\beta_{n,j}|^2=1\eqno(23)
$$
$(j=1,2)$.
Relations (19), (21) and (23) yield
$$
\begin{array}{c}
||\mathbf y_{n,j}||_\mathbb{H}^2=|\alpha_{n,j}|^2\int_0^\pi e^{-2Im\lambda_{n,j}x}dx+|\beta_{n,j}|^2\int_0^\pi e^{2Im\lambda_{n,j}x}dx+o(1)=\\
=|\alpha_{n,j}|^2h_{n,j}+|\beta_{n,j}|^2g_{n,j}+o(1),
\end{array}
$$
where
$$
h_{n,j}=\int_0^\pi e^{-2Im\lambda_{n,j}x}dx,\quad g_{n,j}=\int_0^\pi e^{2Im\lambda_{n,j}x}dx.
$$
This implies the relations
$$
0<c_1<||\mathbf y_{n,j}||_\mathbb{H}<c_2,\eqno(24)
$$

$$
\int_0^\pi y_{n,j}^{[1]}(x)y_{n,j}^{[2]}(x)dx=\pi\alpha_{n,j}\beta_{n,j}+o(1).\eqno(25)
$$
If condition (22) holds, then, it follows from (23) that
$$
|\alpha_{n,1}|,|\beta_{n,1}|>c_3.\eqno(26)
$$
Since
$$
\langle\mathbf{y}_{n,1},\mathbf{z}_{n,2}\rangle=\int_0^\pi(y_{n,1}^{[1]} y_{n,2}^{[2]}+y_{n,1}^{[2]} y_{n,2}^{[1]})dx=0,\eqno(27)
$$
then, by virtue of (27) and (21) we have
$$
\alpha_{n,1}\beta_{n,2}+\alpha_{n,2}\beta_{n,1}=o(1).\eqno(28)
$$
Let us prove that for all sufficiently large $|n|$
$$
|\alpha_{n,2}|>c_4>0.
$$
Indeed, otherwise, there exists a sequence of $n_k$ such that
$$
\lim_{k\to\infty}\alpha_{n_k,2}=0,
$$
this, together with (23), (26) and (28) implies
$$
\lim_{k\to\infty}\beta_{n_k,2}=0,
$$
which contradicts normalization condition (23).
Similarly, for all sufficiently large $|n|$
$$
|\beta_{n,2}|>c_5>0.
$$
These observations together with (25) and (26) imply that for all sufficiently large $|n|$
$$
\biggl|\int_0^\pi y_{n,j}^{[1]}(x) y_{n,j}^{[2]}(x)dx\biggl|>c_6>0.
$$
It follows from the last inequality and (20) that
$$
|\gamma_{n,j}|<c_7.
$$
This and (24) imply
$$
||\mathbf{y}_{n,j}||||\mathbf{z}_{n,j}||=|\gamma_{n,j}|||\mathbf{y}_{n,j}||_\mathbb{H}^2<c_8.
$$
Since the root function system of problem (18) forms a a Riesz basis with parentheses, it follows from the last inequality
that the indicated system is a usual  Riesz basis in $\mathbb{H}$.

If condition (22) does not hold, then, by virtue of (23)
$$
\underline\lim_{|n|\to\infty}|\alpha_{n,1}||\beta_{n,1}|=0.
$$
This and (25) imply
$$
\underline\lim_{|n|\to\infty}\biggl|\int_0^\pi y_{n,1}^{[1]}(x) y_{n,1}^{[2]}(x)dx\biggl|=0
$$
and, according to (20)
$$
\overline\lim_{|n|\to\infty}|\gamma_{n,1}|=\infty.
$$
It follows from the last inequality and (24) that
$$
\overline\lim_{|n|\to\infty}||\mathbf{y}_{n,1}||_\mathbb{H}||\mathbf{z}_{n,1}||_\mathbb{H}=
\overline\lim_{|n|\to\infty}|\gamma_{n,1}|||\mathbf{y}_{n,1}||_\mathbb{H}^2=\infty,
$$
hence, the root function system of problem (18) is not even a usual basis in $\mathbb{H}$.

{\bf Lemma 2.} {\it The root function system of problem (18) forms a Riesz basis if and only if there exist constants $C_1, C_2>0$  such that
$$
C_1<\frac{|e_{12}(\pi,\lambda_{n,1})|}{|e_{21}(\pi,\lambda_{n,1})|}<C_2.
$$
}
Proof. It follows from (2) that the characteristic equation of problem (18) has the form
$$
a^2e_{11}(\pi,\lambda)+2a+e_{22}(\pi,\lambda)=0.\eqno(30)
$$
Estimates (7) imply that in the strip $\Pi$ for all sufficiently large $|\lambda|$ $e_{11}(\pi,\lambda)\ne0$.
Multiplying  equation (30) by $e_{11}(\pi,\lambda)$  we obtain
$$
a^2e^2_{11}(\pi,\lambda)+2ae_{11}(\pi,\lambda)+e_{11}(\pi,\lambda)e_{22}(\pi,\lambda)=0,\eqno(31)
$$
which implies
$$
\biggl(e_{11}(\pi,\lambda)+\frac{1}{a}+\frac{\sqrt{D(\lambda)}}{a}\biggl)\biggl(e_{11}(\pi,\lambda)+\frac{1}{a}-\frac{\sqrt{D(\lambda)}}{a}\biggl)=0,\eqno(32)
$$
where $D(\lambda)=1-e_{11}(\pi,\lambda)e_{22}(\pi,\lambda)=-e_{12}(\pi,\lambda)e_{21}(\pi,\lambda)$. Since $\lambda_{n,1}\ne\lambda_{n,2}$, then
$D(\lambda_{n,1})\ne0$.
Suppose, for example, that the first factor in  (32) has a root $\lambda_{n,1}$. This,  together with (30) implies
$$
e_{11}(\pi,\lambda_{n,1})+\frac{1}{a}=-\frac{\sqrt{D(\lambda_{n,1})}}{a},\quad
e_{22}(\pi,\lambda_{n,1})+a=-a(ae_{11}(\pi,\lambda_{n,1})+1).\eqno(33)
$$
It follows from (33) that
$$
e_{22}(\pi,\lambda_{n,1})+a=a\sqrt{D(\lambda_{n,1})}.\eqno(34)
$$
Denote
$$
 \mathbf{u^j}_{n,1}(x)=\begin{pmatrix}
e_{1j}(x,\lambda_{n,1})\\
e_{2j}(x,\lambda_{n,1})
\end{pmatrix}
$$
$j=1,2$.
Then, according to [18 p. 84], the function

$$
\mathbf{u}_{n,1}(x)=V_2(\mathbf{u^2}_{n,1}(x))\mathbf{u^1}_{n,1}(x)-V_2(\mathbf{u^1}_{n,1}(x))\mathbf{u^2}_{n,1}(x)
$$
is an eigenfunction, corresponding to the eigenvalue $\lambda_{n,1}$. It follows from (33) and (34) that
$$
\begin{array}{c}
\mathbf{u}_{n,1}(x)=((a+e_{22}(\pi,\lambda_{n,1}))\mathbf{u^1}_{n,1}(x)-e_{21}(\pi,\lambda_{n,1})\mathbf{u^2}_{n,1}(x)=\\=
a\sqrt{D(\lambda_{n,1}})\mathbf{u^1}_{n,1}(x)-e_{21}(\pi,\lambda_{n,1})\mathbf{u^2}_{n,1}(x)=\\=
a\sqrt{-e_{12}(\pi,\lambda_{n,1})e_{21}(\pi,\lambda_{n,1}})\mathbf{u^1}_{n,1}(x)-e_{21}(\pi,\lambda_{n,1})\mathbf{u^2}_{n,1}(x).
\end{array}
$$
The assertion of Lemma 2 follows from the last equality and Lemma 1.

Let us find the asymptotic behavior  of the functions $e_{12}(\pi,\lambda),e_{21}(\pi,\lambda)$ if $|\lambda|\to\infty$, $\lambda\in\Pi$.
In [19] the Dirac system was considered in the form
$$
 \tilde L\mathbf{y}=\tilde B\mathbf{y}'+\tilde P\mathbf{y},
 $$
 where
 $$
 \tilde B=\begin{pmatrix}
 0&-1\\
 -1&0
 \end{pmatrix},\quad \tilde P=\begin{pmatrix}
 p(x)&r(x)\\
 r(x)&-p(x)
 \end{pmatrix}.
$$
In particular, it was established that if $p(x), r(x)\in W_2^1[0,\pi]$, then elements $y_{ij}$ of the fundamental solution $y(\lambda,x)$
 of the equation
 $$
 \tilde L\mathbf{y}=\lambda \mathbf{y} \eqno(35)
 $$
 with boundary conditions $y(\lambda,0)=I$
can be represented in the form $y_{ij}=\hat y_{ij}/w(\lambda)$, where
$$
\begin{array}{c}
\hat y_{11}=\elx u^+(\lambda,x)[1+\sigma^-(\lambda,x)][1+\sigma^+(-\lambda,0)]+\\+
\elxx u^-(-\lambda,x)[1-\sigma^+(-\lambda,x)][1-\sigma^-(\lambda,0)],
\end{array}\eqno(36)
$$
$$
\begin{array}{c}
i\hat y_{12}=-\elx u^+(\lambda,x)[1-\sigma^-(\lambda,x)][1+\sigma^+(-\lambda,0)]+\\+
\elxx u^-(-\lambda,x)[1+\sigma^+(-\lambda,x)][1-\sigma^-(\lambda,0)],
\end{array}\eqno(37)
$$
$$
\begin{array}{c}
i\hat y_{21}=\elx u^+(\lambda,x)[1+\sigma^-(\lambda,x)][1-\sigma^+(-\lambda,0)]-\\-
\elxx u^-(-\lambda,x)[1-\sigma^+(-\lambda,x)][1+\sigma^-(\lambda,0)],
\end{array}\eqno(38)
$$
$$
\begin{array}{c}
\hat y_{22}=\elx u^+(\lambda,x)[1-\sigma^-(\lambda,x)][1-\sigma^+(-\lambda,0)]+\\+
\elxx u^-(-\lambda,x)[1+\sigma^+(-\lambda,x)][1+\sigma^-(\lambda,0)],
\end{array}\eqno(39)
$$
$$
w(\lambda)=2(1+\sigma^+(-\lambda,0)\sigma^-(\lambda,0)),\eqno(40)
$$
where
$$
u^{\pm}(\lambda,x)=1+\frac{b_1^\pm(x)}{2i\lambda}+O(\lambda^{-2}),\eqno(41)
$$
$$
b_1^\pm(x)=\int_0^x(p^2(x)+r^2(x))dx,\eqno(42)
$$
$$
\sigma^\pm(\lambda,x)=1+\frac{\sigma_1^\pm(x)}{2i\lambda}+\frac{\tilde\sigma^\pm(\lambda,x)}{2i\lambda}+O(\lambda^{-2}),\eqno(43)
$$
$$
\sigma_1^\pm(x)=ip(x)\mp r(x),\quad \tilde\sigma^\pm(\lambda,x)=\int_0^x\sigma_2^\pm(x-t)e^{-2i\lambda t}dt,\quad\sigma_2^\pm(x)=\pm r'(x)-ip'(x).\eqno(44)
$$

Let
$\mathbf{y}=(y_1,y_2)$ be a solution of problem (35). We replace $u_1=y_1+iy_2, u_2=y_1-iy_2$. Simple computations show that
the vector $\mathbf{u}=(u_1,u_2)$ is a solution of the system $L\mathbf{u}=\lambda\mathbf{u}$, where $P(x)=p(x)+ir(x), Q(x)=p(x)-ir(x)$.
Consider the matrix
$$
Y(\lambda,x)=\begin{pmatrix}
y_{11}+iy_{21}&y_{12}+iy_{22}\\
y_{11}-iy_{21} &y_{12}-iy_{22}
 \end{pmatrix}.
$$
Obviously,
$$
Y(\lambda,0)=\begin{pmatrix}
 1&i\\
1 &-i
 \end{pmatrix},
$$
which implies that
$$
Y^{-1}(\lambda,0)=\begin{pmatrix}
 1&1\\
-i &i
 \end{pmatrix}/2.
$$
It follows from this and [19, p. 66] that
$$
\begin{array}{c}
E(x,\lambda)=\begin{pmatrix}
e_{11}(x,\lambda)&e_{12}(x,\lambda)\\
e_{21}(x,\lambda)&e_{22}(x,\lambda)
\end{pmatrix}=Y(\lambda,x)Y^{-1}(\lambda,0)=\\
=\frac{1}{2}\begin{pmatrix}
y_{11}+y_{22}+i(y_{21}-y_{12})&y_{11}-y_{22}+i(y_{12}+y_{21})\\
y_{11}-y_{22}-i(y_{12}+y_{21}) &y_{11}+y_{22}-i(y_{21}-y_{12})
 \end{pmatrix}.
 \end{array}\eqno(45)
$$
Relations (36-40) imply that
$$
\begin{array}{c}
e_{12}(x,\lambda)=\frac{1}{2}(y_{11}-y_{22}+i(y_{12}+y_{21}))=\\
=\frac{1}{2w(\lambda)}(\elx u^+(\lambda,x)[1+\sigma^-(\lambda,x)][1+\sigma^+(-\lambda,0)]+\\+
\elxx u^-(-\lambda,x)[1-\sigma^+(-\lambda,x)][1-\sigma^-(\lambda,0)]-\\-
(\elx u^+(\lambda,x)[1-\sigma^-(\lambda,x)][1-\sigma^+(-\lambda,0)]+\\+
\elxx u^-(-\lambda,x)[1+\sigma^+(-\lambda,x)][1+\sigma^-(\lambda,0)])+\\+
(-\elx u^+(\lambda,x)[1-\sigma^-(\lambda,x)][1+\sigma^+(-\lambda,0)]+\\+
\elxx u^-(-\lambda,x)[1+\sigma^+(-\lambda,x)][1-\sigma^-(\lambda,0)])+\\+
\elx u^+(\lambda,x)[1+\sigma^-(\lambda,x)][1-\sigma^+(-\lambda,0)]-\\-
\elxx u^-(-\lambda,x)[1-\sigma^+(-\lambda,x)][1+\sigma^-(\lambda,0)],
\end{array}\eqno(46)
$$
$$
\begin{array}{c}
e_{21}(x,\lambda)=\frac{1}{2}(y_{11}-y_{22}-i(y_{12}+y_{21}))=\\
=\frac{1}{2w(\lambda)}(\elx u^+(\lambda,x)[1+\sigma^-(\lambda,x)][1+\sigma^+(-\lambda,0)]+\\+
\elxx u^-(-\lambda,x)[1-\sigma^+(-\lambda,x)][1-\sigma^-(\lambda,0)]-\\-
(\elx u^+(\lambda,x)[1-\sigma^-(\lambda,x)][1-\sigma^+(-\lambda,0)]+\\+
\elxx u^-(-\lambda,x)[1+\sigma^+(-\lambda,x)][1+\sigma^-(\lambda,0)])-\\-
[(-\elx u^+(\lambda,x)[1-\sigma^-(\lambda,x)][1+\sigma^+(-\lambda,0)]+\\+
\elxx u^-(-\lambda,x)[1+\sigma^+(-\lambda,x)][1-\sigma^-(\lambda,0)])+\\+
\elx u^+(\lambda,x)[1+\sigma^-(\lambda,x)][1-\sigma^+(-\lambda,0)]-\\-
\elxx u^-(-\lambda,x)[1-\sigma^+(-\lambda,x)][1+\sigma^-(\lambda,0)]].
\end{array}\eqno(47)
$$
Substituting in equalities (46) and (47) asymptotic representations (41-44), we obtain
$$
\begin{array}{c}
e_{12}(\pi,\lambda)=\frac{1}{i\lambda w(\lambda)}(e^{i\pi\lambda}(\sigma_1^-(\pi)+\tilde\sigma^-(\lambda,\pi)+O(\lambda^{-1}))-e^{-i\pi\lambda}(\sigma_1^-(0)+O(\lambda^{-1}))=\\
=\frac{1}{\lambda w(\lambda)}[e^{i\pi\lambda}(Q(\pi)+e^{-2i\lambda}\int_0^\pi Q'(t)e^{2i\lambda}dt+O(\lambda^{-1}))-e^{-i\pi\lambda}(Q(0)+O(\lambda^{-1})].
\end{array}\eqno(48)
$$
$$
\begin{array}{c}
e_{21}(\pi,\lambda)=\frac{1}{i\lambda w(\lambda)}(-e^{i\pi\lambda}(\sigma_1^+(0)+O(\lambda^{-1}))+e^{-i\pi\lambda}(\sigma_1^+(\pi)+\tilde\sigma^+(-\lambda,\pi)+O(\lambda^{-1}))=\\
=\frac{1}{\lambda w(\lambda)}[-e^{i\pi\lambda}(P(0)+O(\lambda^{-1}))+e^{-i\pi\lambda}(P(\pi)+e^{2i\lambda}\int_0^\pi P'(t)e^{-2i\lambda t}dt+
O(\lambda^{-1}))].
\end{array}\eqno(49)
$$
Denote by $\Psi$ the set of pair of functions $(P(x),Q(x))\in L_1(0,\pi)\oplus L_1(0,\pi)$ such that the root function system of problem (18) forms a Riesz basis in $\mathbb{H}$, $\overline{\Psi}=(L_1(0,\pi)\oplus L_1(0,\pi))\setminus\Psi$,

{\bf Theorem 2.} {\it The sets $\Psi$ and
 $\overline{\Psi}$ are everywhere dense in $L_1(0,\pi)\oplus L_1(0,\pi)$.}

 Proof. Suppose functions $f_i(x)\in L_1(0,\pi)$, $(i=1,2)$, $\varepsilon>0$. It is readily seen that there exist functions $\hat P(x), \hat Q(x)$ such that $\hat P(x)\in C^1[0,\pi]$, $\hat P(0)=0$, $\hat P(\pi)\ne0$, $||f_1(x)-\hat P(x)||_{L_1(0,\pi)}<\varepsilon$/2, $\hat Q(x) \in C^1[0,\pi]$, $\hat Q(0)=0$, $\hat Q(\pi)\ne0$, $||f_2(x)-\hat Q(x)||_{L_1(0,\pi)}<\varepsilon$/2. Since for any $n\in\mathbb{Z}$ the inequality $0<c_1<|e^{i\pi\lambda_{n,j}}|<c_2$ takes place, then it follows from (48), (49) that the elements of the fundamental matrix of problem (18) with the potential  $\hat P(x), \hat Q(x)$ satisfy the inequalities
 $$
\frac{c_3}{|\lambda_{n,j}|}<|e_{12}(\pi,\lambda_{n,j})|<\frac{c_4}{|\lambda_{n,j}|},\quad
\frac{c_3}{|\lambda_{n,j}|}<|e_{21}(\pi,\lambda_{n,j})|<\frac{c_4}{|\lambda_{n,j}|}.
$$
The last inequality and Lemma 2 imply that the set $\Psi$ is everywhere dense in $L_1(0,\pi)\oplus L_1(0,\pi)$.

Obviously, there exists a function $S(x)\in C^1[0,\pi]$ such that $S(0)=S(\pi)=S'(0)=S'(\pi)=0$, $||f_1(x)-S(x)||_{L_1(0,\pi)}<\varepsilon/10$.
Consider the Fourier-series expansion of the function $S(x)$ on the segment $[0,\pi]$
$$
S(x)=\alpha_0+\sum_{m=1}^\infty\alpha_m e^{2imx}+\beta_me^{-2imx}.
$$
It is well known that there exists a number $N$ such that $S(x)=S_N(x)+R_N(x)$, where
$$S_N(x)=\alpha_0+\sum_{m=1}^N\alpha_m e^{2imx}+\beta_me^{-2imx},
$$
and $|R_N(x)|<\varepsilon/10$. Evidently, $|S_N(0)|<\varepsilon/10$, $|S_N(\pi)|<\varepsilon/10$,
$$
S_N'(x)=2i\sum_{m=1}^Nm(\alpha_m e^{2imx}-\beta_me^{-2imx}).
$$
Let $a_0=N$, $a_{k+1}=Ca_k$, $k=0,1,\ldots$, where $C=[e^{2|w_0|\pi}+100]N^2[1/\varepsilon^2+1]$.
It is readily seen that
$$
\sum_{j=1, j\ne k}^\infty\frac{1}{|a_j-a_k|}=\sum_{j=1}^{k-1}\frac{1}{a_k-a_j}+\sum_{j=k+1}^\infty\frac{1}{a_j-a_k}
<c_1\biggl(\frac{k}{a_k}+\frac{1}{Ca_k}\biggl)<\frac{c_2}{a_k^{2/3}}.\eqno(50)
$$
Denote
$$
\theta(x)=e^{2iw_0 x}\sum_{k=1}^\infty\frac{e^{4ia_kx}}{\sqrt{a_k}}.
$$
It is easily shown that
$$
|\theta(x)|<\varepsilon/(10\pi).\eqno(51)
$$

Let us estimate the integral
$$
I_0=\int_0^\pi\theta(t)e^{-2i\lambda_{a_k,1} t}dt=I_1+I_2,\eqno(52)
$$
where
$$
I_1=\frac{1}{\sqrt{a_k}}\int_0^\pi e^{-2i\varepsilon_{a_k,1}t}dt,\quad
I_2=\sum_{j=1, j\ne k}^\infty\frac{1}{\sqrt{a_j}}
\int_0^\pi e^{2i(2a_j-2a_k-\varepsilon_{a_k,1})t}dt.
$$
It is obvious that
$$
\frac{3}{\sqrt{a_k}}<|I_1|<\frac{4}{\sqrt{a_k}}.\eqno(53)
$$
It follows from (50) that
$$
|I_2|<c_3\sum_{j=1, j\ne k}^\infty\frac{1}{\sqrt{a_j}}\frac{1}{|a_j-a_k|}<\frac{c_4}{a_k^{2/3}}\eqno(54)
$$
and
$$
\biggl|\int_0^\pi S_N'(t)e^{-2i\lambda_{a_k,1} t}dt\biggl|
<\frac{c_5N^2}{a_k-N}<\frac{c_6}{a_k^{2/3}}.\eqno(56)
$$
Estimates (52-54) imply that the following inequality holds for all sufficiently large $k$
$$
\frac{2}{\sqrt{a_k}}<|I_1|<\frac{5}{\sqrt{a_k}}.\eqno(55)
$$

Denote
$$
F(x)=S_N(x)+\int_0^x\theta(t)dt,\quad F_0(x)=\frac{F(\pi)-F(0)}{\pi}x+F(0),\quad\tilde P(x)=F(x)-F_0(x).
$$
Evidently,
$\tilde P(0)=\tilde P(\pi)=0$.
It follows from (51) that $|F(0)|<\varepsilon/10, \quad|F(\pi)|<\varepsilon/5$. These estimates imply
$$
\biggl|\int_0^\pi F_0'(t)e^{-2i\lambda_{a_k,1} t}dt\biggl|<\frac{c_7}{a_k}.
$$
This inequality, together with (55) and  (56) gives the inequality
$$
\frac{c_8}{\sqrt{a_k}}<\biggl|\int_0^\pi \tilde P'(t)e^{-2i\lambda_{a_k,1} t}dt\biggl|<\frac{c_9}{\sqrt{a_k}}.
$$
The last relation and (49) imply that
$$
\frac{c_{10}}{a_k\sqrt{a_k}}<|e_{21}(\pi,\lambda_{a_k,1})|<\frac{c_{11}}{a_k\sqrt{a_k}}.
$$
Thus, according to Lemma 2 the necessary condition of basis property fails  for problem (18) with the potential $(\tilde P(x), \hat Q(x))$.
 Obviously,
$$
||S(x)-\tilde P(x)||_{L_1(0,\pi)}<2\varepsilon/5,
$$
therefore,
$$
||f_1(x)-\tilde P(x)||_{L_1(0,\pi)}<\varepsilon/2.
$$
Theorem 2 is proved.

\medskip
\medskip
\medskip

\centerline {REFERENCES}
\medskip
\medskip
\medskip

\begin{itemize}
\item[1.] A.M. Savchuk, A.A. Shkalikov. The Dirac Operator with Complex-Valued Summable Potential. Math. Notes. 2014. V. 96. No. 5.
P. 777-810.

\item[2.] A.M. Savchuk, I. V. Sadovnichaya. Riesz basis property with parentheses for the Dirac system with summable potential. Sovrem. Math. Fundam. Napravl. 2015. V. 58. P. 128-152.

\item[3.] A. Lunyov, M. Malamud. On the Riesz basis property of root vectors system for $2\times2$ Dirac type operators. J. Math. Anal. Appl.,
2016, V. 441. No. 1. P. 57-103.

\item[4.] A. Lunyov, M. Malamud. On the Riesz basis property of root vectors system for $2\times2$ Dirac type systems. Dokl. Math. 2014. V. 90. No. 2. P. 556-561.

\item[5.] P. Djakov,  B. Mityagin. Unconditional Convergence of Spectral Decompositions of 1D Dirac operators with Regular Boundary Conditions. Indiana Univ. Math. J. 2012.  V. 61. No. 1. P. 359-398.

\item[6.] I. Ts. Gohberg and M.G. Krein. Introduction to the Theory of Linear Non-self-adjoint Operators in Hilbert Space. (Am. Math, Soc., Providence, RI, 1969, M., Mir, 1967).

\item[7.] M.V. Keldysh. On the completeness of the eigenfunctions of some classes of non-self-adjoint linear operators. Russ. Math. Surv. 1971. V. 26. No. 4. P. 15-41.

\item[8.] K. Iosida. Functional Analysis (Springer-Verlag, Berlin, 1965, Mir, Moscow, 1967).

\item[9.] P. Djakov,  B. Mityagin. Instability zones of one-dimensional periodic Schrodinger and Dirac operators. Uspekhi Mat. Nauk 61 (2006) No. 4, 77-182.; English transl. Russian Math. Surveys 61 (2006) No 4, 663-766.

\item[10.] P. Djakov,  B. Mityagin. Bari-Markus property for Riesz projections of 1D periodic Dirac operators, Math. Nachr.  2010. V. 283. No. 3. P. 443-462.

\item[10.] P. Djakov,  B. Mityagin. 1D Dirac operators with special periodic potentials. Bull. Pol. Acad. Sci. Math. 2012. V. 60.  No. 1. P. 59-75.

\item[11.] P. Djakov,  B. Mityagin. Equiconvergence of spectral decompositions of 1D Dirac operators with regular boundary conditions.
J. Approx. Theory. 2012. V. 164.  No. 7. P. 879-927.

\item[11.] P. Djakov,  B. Mityagin. Criteria for existence of Riesz bases consisting of root functions of Hill and 1D Dirac operators.
J. Funct. Anal. 2012. V. 263.  No. 8. P. 2300-2332.

\item[11.] P. Djakov,  B. Mityagin. Riesz bases consisting of root functions of 1D Dirac operators. Proc. AMS. 2013. V. 141. No. 4. P. 1361-1375.

\item[12.] B. Mityagin. Spectral Expansions of One-dimensional Periodic Dirac Operators. Dynamics of PDE. 2004. V. 1. No 2. P. 125-191.

\item[13.] B. Mityagin. Convergence of expansions in eigenfunctions of the Dirac operator. Dokl. Math. 2003. V. 68. No 3. P. 388-391, Trans. from Dokl. Akad. Nauk. 2003. V. 393. No 4. P. 456-459.

\item[14.] A.G. Baskakov, A.V. Derbyshev, A.O. Shcherbakov. The method of similar operators in the spectral analysis of non-self-adjoint Dirac operators with non-smooth potentials. Izvestiya Mathematics. 2011. V. 75. No. 3. P. 3-28.

\item[15.] M.Sh. Burlutskaya, V.V. Kornev, A.P. Khromov. Dirac system with nondifferentiable potential and periodic boundary conditions.
      Zh. Vychisl. Mat. Mat. Fiz. 2012. V. 52. No. 8. P. 1621-1632.

\item[16.] V.V. Kornev, A.P. Khromov. Dirac system with nondifferentiable potential and antiperiodic boundary conditions. Izv. Sarat. Mat. Mekx. Inform. 2013. V. 13. No. 3. P. 28-35.

\item[17.] I. Arslan. Characterization of the potential smoothness of one-dimensional Dirac operator subject to general boundary conditions and its basis property. J. Math. Anal. Appl. 2017. V. 447. No. 1. P. 84-108.

\item[18.] M. A. Naimark. Linear Differential Operators. (F. Ungar Publ. New York, 1967, Nauka, Moscow, 1969).

\item[19.] V. A. Marchenko. Sturm-Liouville Operators and Their Applications. ( Kiev, 1977) [in Russian]).

\end{itemize}

\medskip
\medskip
\medskip
\medskip
\medskip
\medskip
E-mail: alexmakin@yandex.ru

}
\end{document}